\documentstyle[12pt]{article}

\makeatletter
\newcommand{\newoperator}[2]{\@ifdefinable#1{\def#1{\mathop{#2}\nolimits}}}
\makeatother

\newcommand\Differential[2]{\bigD#1\ #2}

\newcommand\Dom[1]{\mbox{Dom}\left(#1\right)}

\newcommand\Emu[2]{\expectation_{#1}\left[#2\right]}

\newcommand\integers{\IN}

\newcommand\nDifferential[3]{\bigD^{#1}#2\ #3}

\newcommand\norm[1]{\Vert#1\Vert}

\newcommand\reals{\IR}
\newcommand\realvectors[1]{\IR^#1}

\newcommand\set[1]{\left\{#1\right\}}

\newoperator{\AND}{\rm AND}
\newoperator{\area}{\rm area}
\newoperator{\bigD}{\rm D}
\newoperator{\card}{\rm card}
\newoperator{\cov}{\rm Cov}
\newoperator{\determinant}{\rm det}
\newoperator{\expectation}{\rm \IE}
\newoperator{\hess}{\rm Hess}
\newoperator{\NOR}{\rm NOR}
\newoperator{\NOT}{\rm NOT}
\newoperator{\OR}{\rm OR}
\newoperator{\parti}{\Cal P}
\newoperator{\probability}{\IP}
\newoperator{\smallO}{\rm o}
\newoperator{\var}{\rm Var}
\newoperator{\XOR}{\rm XOR}
\newtheorem{theorem}{Theorem}
\newtheorem{definition}[theorem]{Definition}
\newtheorem{proposition}[theorem]{Proposition}

\newenvironment{example}
{\refstepcounter{theorem}\par\noindent{\em Example \thetheorem.}\
}{\par}

\newcommand\Cal[1]{{\cal #1}}
\newcommand\emmex{\Cal M(X, \Cal X, \mu)}

\newcommand\formula[1]{(\ref{eq:#1})}
\newcommand\leftstar{{^*}\!}

\newcommand\proof{\par\noindent{\em Proof.\ }}
\newcommand\QED{$\Box$\par\bigskip}

\newcommand\IN{I\!\!N}
\newcommand\IR{I\!\!R}

\newcommand\IP{I\!\!P}
\newcommand\IE{I\!\!E}

\renewcommand{\span}{{\rm span}}
\newcommand{\inverse}[1]{{\textstyle\frac{1}{#1}}}
\newcommand{\half}{\inverse{2}}

\title{{\normalsize (Preprint 4/1996, Dipartimento di Matematica, Universit\`a di Padova, 1996)}\\ Projecting the Fokker-Planck Equation onto
a finite dimensional exponential family }

\author{Damiano Brigo \thanks{While working on this article
the author was supported by a senior fellowship
of the "Istituto Nazionale di Alta Matematica F. Severi" (Rome, Italy) }
\\
Dipartimento di Matematica  \\
Universit\`a di Padova
\\
Via Belzoni 7\\
35131 Padova\\
Italy\\
brigo@pdmat1.unipd.it
\and
Giovanni Pistone\\
Dipartimento di Matematica\\ Politecnico di Torino\\ Corso Duca
degli Abruzzi 24\\ 10129 Torino\\ Italy\\ pistone@polito.it}
\date{}
\begin{document}

\maketitle

\thispagestyle{empty}

\begin{abstract}
In the present paper we discuss problems concerning
evolutions of densities related to It\^o diffusions in the framework
of the statistical exponential manifold. We develop a rigorous
approach to the problem, and we particularize it to
the orthogonal projection of the evolution of
the density of a diffusion
process onto a finite dimensional exponential manifold.
It has been shown by D. Brigo (1996) that the projected evolution
 can always be interpreted as the evolution of the density
of a different diffusion process.
We give also a compactness result when the dimension of the
exponential family increases, as a first step towards a convergence
result
to be investigated in the future.
The infinite dimensional
exponential manifold structure introduced by G. Pistone and C. Sempi
is used and some examples are given.
\end{abstract}


\section{Introduction}
This paper moves both from the differential geometric approach to
nonlinear filtering as developed by Brigo, Hanzon and LeGland
\cite{BrHaLe} and from the rigorous approach to the construction
of a differential geometric structure in the infinite dimensional
space of probability measures given in Pistone and Sempi \cite{PSS92},
see also Pistone and Rogantin \cite{PSR96}.
The solution of the filtering problem is a stochastic PDE
which can be seen as a generalization of the Fokker--Planck
equation (FPE) expressing the density of a diffusion process.
This filtering equation is called the Kushner--Stratonovich equation.
In \cite{BrHaLe} the Fisher metric is used to
project the Kushner--Stratonovich equation
onto a finite dimensional exponential manifold of probability densities.
This method can be used also for the simpler FPE.
In the present paper we discuss the geometric approach to
problems concerning finite dimensionality of densities related to
stochastic differential equations given by It\^o diffusions.
Part of these results were already given in \cite{Brigo4},
and we shortly present them here in the framework of
Pistone and Sempi \cite{PSS92}.
This approach is different from the one
adopted for example in Brigo, Hanzon and LeGland \cite{BrHaLe}
or in Brigo \cite{Brigo4},
since it uses the exponential manifold structure rather than
the $L_2$ derivation.
The $L_2$ structure is obtained by mapping densities into their
square roots. We show that this map yields a regular
$C^\infty$ parametrization but it is not a chart for the
infinite dimensional manifold of densities.
%
%For an extensive treatment of the extension of the Fisher metric to an
%infinite dimensional family of densities see \cite{PSS92}.
In the present paper we consider the projection in Fisher metric
of the density--evolution of a diffusion process onto an exponential
manifold.
Such projection is obtained via the projected FPE. We examine the
projected density--evolution and discuss problems related to finite
dimensionality,
giving some examples.
We recall from Brigo \cite{Brigo4} that the projected density--evolution
can always be
interpreted as the density--evolution
of a different diffusion process. We conclude by
giving a first step for future investigations on the following
convergence problem:
is it possible to prove that
the projected density converges to the original one when the
dimension of the exponential manifold on which we project tends
to infinity?

\section{The exponential statistical manifold of positive
probability densities}
In the present section we give a summary of the construction of the
non-parametric exponential statistical manifold as developed in
\cite{PSS92} and \cite{PSR96}.
In those papers it is shown that the
definition of statistical manifold as introduced by Dawid, Efron, Amari
and
others, and systematically presented in \cite{MRR93},
can be given in a non parametric setting using the
framework of the theory of manifolds modeled on Banach spaces, as
introduced for example in Lang, \cite{LNG95}.

We consider a measure space $(X, \Cal X, \mu)$, where $\mu$ is a
reference measure, and the set $\emmex$ of the a.s. strictly
positive densities w. r. t. some measure equivalent to $\mu$.  We define
on the set $\emmex$ a topology such that $\emmex$ is an Hausdorff space
(i.e. points can be
separated by open sets). Then we shall construct a covering of $\emmex$
with open
sets $\Cal U_p$, $p\in\Cal U_p$, $p\in\emmex$, and a corresponding
family of Banach spaces $B_p$, with norms $||\cdot||_p$, $p\in
\emmex$, such that each density $q\in\Cal U_p$ is represented by a
coordinate $s_p(q) \in B_p$.

We shall use the notations

\begin{eqnarray} s_p&:& \Cal U_p \to \Cal V_p \subset B_p
\label{eq:chart} \\ e_p&:& \Cal V_p \to \Cal U_p \subset \emmex
\label{eq:patch} \end{eqnarray} to denote respectively the {\em
charts}, i.e. the mappings from points to coordinates, and the {\em
patches}, i.e. the mappings from coordinates to points.

Following the use in differential geometry, we say that $\{(\Cal U_p,
s_p) \ :
\ p\in\emmex \}$ is an {\em atlas} if all the space
is covered by charts;  if moreover each of the {\em change of
coordinates}
$$s_{p_2}\circ e_{p_1} : s_{p_1}\left(\Cal U_{p_1} \cap \Cal
U_{p_2}\right) \to e_{p_2}\left(\Cal U_{p_1} \cap \Cal
U_{p_2}\right)$$ is a diffeomorphism of some regularity
between open sets, the atlas has that
regularity. In such a case the atlas, augmented with all the
compatible charts, defines the manifold, see Lang \cite{LNG95}.

In our case we shall introduce a very special manifold, such that
the change of coordinates are actually affine functions ---i.e. they
differ from a linear function by a constant---, but we will keep a
weaker regularity, namely the $C^\infty$--regularity
(differentiability of any order) for compatible charts.

We shall denote by $\Emu{p\cdot\mu}{\cdot}$ the expectation w.r.t.
the probability measure $p\cdot\mu$ (where $p\cdot\mu(dx)=p(x)\mu(dx)$);
if
there
is no ambiguity we will use the notation $\Emu p {\cdot}$.

First we define the topology as follows. For simplicity we give only
the definition of convergence of sequences.
The sequence
$(p_n)_{n\in\integers}$ in $\emmex$ is $e$--convergent
(exponentially convergent) to $p$ if $(p_n)_{n\in\integers}$ tends
to $p$ in $\mu$--probability as $n\to\infty$ and moreover the
sequences $(p_n/p)_{n\in\integers}$ and $(p/p_n)_{n\in\integers}$
are eventually bounded in each $L^\alpha(p)$, $\alpha >1$, i.e.
\begin{displaymath}\forall \alpha > 1\quad \limsup_{n\to\infty}
E_p\left[\left({p_n\over p}\right)^\alpha\right]<+\infty,
\quad\limsup_{n\to\infty} E_p\left[\left({p\over
p_n}\right)^\alpha\right]<+\infty.\end{displaymath}

Now we shall introduce the Banach spaces on which the statistical
manifold is modeled. We give a definition that shows how they are
connected with well-known statistical objects.
For each density $p \in \emmex$,
the {\em Cramer class at $p$} is the set of all random
variables $u$ on $X$ such that the {\em moment generating function}
\begin{displaymath}
\hat{u}_{p}(t) = \int e^{tu} p\,d\mu =
\Emu{p}{e^{tu}},\quad t \in \reals
\end{displaymath}
is finite in a neighborhood of the origin 0. If moreover the expectation
of
$u$ is
zero (the previous condition implies the existence of a finite
expectation), then we shall call the set the {\em centered Cramer class
at $p$}.

The centered
Cramer class at $p$ is a
vector space, and it shall be denoted by $B_p$, i.e.
\begin{displaymath} B_p =\{u \in L^1(p\cdot\mu):0 \in Dom(\hat
u_p)^\circ, \Emu{p}{u}=0\}. \end{displaymath}
It is a Banach space for the norm defined by: \begin{equation} \Vert
u \Vert_p = \inf \left\{r : \Emu{p}{\cosh\left(\frac u r\right) - 1}
\le 1 \right\} \end{equation}

In the previous formula the function $x \mapsto \cosh(x) - 1$ is
a convex function that plays in the theory of the spaces $B_p$ the
same role as the function $x \mapsto |x|^\alpha / \alpha$ in the
theory of Lebesgue spaces $L^\alpha$, $\alpha > 1$. We cite
\cite{KRR58} and \cite{RAR91} as general references.

We will denote by $\leftstar
B_p$ the Banach space of centered random variables of the so called
$x\log x$-class. A random variable $u$ belongs to the $x\log
x$-class $^*B_p$ if and only if it is centered and $(1+u)\log(1+u)$
is $(p\cdot\mu)$-integrable.

Now we give some details about the Banach spaces $B_p$ and
$\leftstar B_p$ which will be useful in the construction of the
statistical manifold.

\begin{proposition}\begin{enumerate}

\item
The dual space of the Banach space $\leftstar B_p$ is isomorphic to
$B_p$,
i.e. if $T$ is a continuous linear operator on $\leftstar B_p$ then
there exists a unique $u\in B_p$ such that $T(k) = \Emu p {ku}$, $k
\in \leftstar B_p$;
that is $(\leftstar B_p)^{\star} \ni T \leftrightarrow u \in B$.

\item All the elements $k$ in $\leftstar B_p$ are identified with an
element of the dual space $B^*_p$ of $B_p$ with the identification
$S(u) = \Emu p {ku}$, but $\leftstar B_p$ is strictly smaller than
$B^*_p$ unless the sample space has a finite number of atoms.

\item Denoting with sub-0 the spaces of centered random variables,
the following continuous inclusions hold true: \begin{displaymath}
L^{\infty}_0(p\cdot\mu) \subset B_p \subset \bigcap_{\alpha>1}
L^\alpha_0(p\cdot\mu) \subset L^\alpha_0(p\cdot\mu) \subset
\leftstar B_p \subset B^*_p. \end{displaymath}
\end{enumerate}
\end{proposition}
The patches of the atlas will be defined on the open ball of radius
1: \begin{displaymath} \Cal V_p = \left\{u \in B_p :
\left\Vert{u}\right\Vert_p < 1\right\} ;\end{displaymath} remark
that the condition $\norm u _p < 1$ is equivalent to the existence
of an $\alpha > 1$ such that $\Emu p {\cosh(\alpha u) - 1} \le 1$,
which in turn implies $\Emu p {e^u} < 4$, see Prop. \ref{prop:kappa}
below.

The {\em moment
generating functional} $G_p: L^{(\cosh\cdot-1)}(p\cdot\mu) \to
{\overline{ R}}_{+} = [0,+\infty]$ is defined by \begin{displaymath}
G_p(u) = \Emu p {e^u} \quad . \end{displaymath}
%
%\begin{proposition}[Properties of the MGF]
%
% The moment generating functional $G_p$
%
%\begin{enumerate}
%
%\item takes value 1 at 0, otherwise is strictly greater than 1, is
%convex and its proper domain $\Dom{G_p} = \{u \in
%L^{(\cosh\cdot-1)}(p\cdot\mu): G_p(u)<\infty \}$ is a convex set
%which contains the open unit ball $\left\{u:\norm u _p < 1\right\}$;
%
%\item is bounded and infinitely Fr\'echet--differentiable on the
%open unit ball $$\left\{u:\norm u _p < 1\right\}$$ with
%differential: \begin{displaymath} \nDifferential n {G_p(u)}
%{(v_1,\dots,v_n)} = \Emu{q}{v_1 \cdots v_n e^{u}}. \end{displaymath}
%
%\end{enumerate} \label{prop:MGF} \end{proposition}
%
The {\em
cumulant generating functional} $K_p : B_p \to [0, +\infty]$ is
defined by
\begin{displaymath} K_p(u) = \log G_p(u). \end{displaymath}
%We remark that we restrict the cumulant generating functional to be
%defined on {\em centered} random variables of the Cramer class at
%$p$.
%
\begin{proposition}[Properties of the CGF] \label{prop:CGF} The
cumulant generating functional $K_p$ has proper domain
$\Dom{G_p}\cap B_p$. If $\Cal V_p$ denotes a subset of the proper domain
then
$K_p$ satisfies the following properties
\begin{enumerate}
\item $K_p$ is 0 at 0, otherwise is strictly positive; is convex and
infinitely Fr\'echet differentiable on $\Cal V_p$. The value at 0 of
the differential of order $n$ is the value of the $n$-th cumulant
under $p$ of the random variable $u$.
\item $\forall u \in \Cal V_p$, $q = e^{u-K_p(u)} \cdot p$ is a
probability density in $\emmex$ and the value of the $n$-th
differential at $u$ in the direction $v$ of $K_p$ is the $n$-th
cumulant of $v$ under $q$:
 \begin{displaymath} \nDifferential n {K_p\left(u\right)} {v^n}
=\left. \frac{d^n}{dt^n} \log \Emu q {e^{tv}}\right|_{t=0}.
\end{displaymath}
\item \label{item:firsttwo} In particular $\frac q p -1 \in
\leftstar B_p$ and
\begin{equation} \begin{array} {c} \Differential {K_p(u)} v = \Emu q
v = \Emu p {\left( \frac q p - 1\right) v} \\ \nDifferential 2
{K_p(u)}{v_1 v_2} = \Emu q {v_1v_2} - \Emu q {v_1} \Emu q {v_2}
\end{array}
\quad .
\end{equation}
\end{enumerate}
\label{prop:kappa}\end{proposition}

Using the definitions introduced so far, it is possible to give a
definition of the non-parametric exponential model as follows.
For each $p$ in $\emmex$
{\em the maximal exponential model} at $p$ is the statistical model
\begin{displaymath} \Cal E_p = \left\{e^{u - K_p(u)} \cdot p: u \in
\Dom{K_p}^\circ, \Emu p u = 0 \right\} \quad . \end{displaymath} The
function \begin{displaymath} B_p \supset \Dom{K_p}^\circ \ni u
\mapsto e^{u - K_p(u)}\cdot p \in \emmex \end{displaymath} is the
likelihood function when the `model parameter' is $u$.

We now have all the elements for the
definition of the atlas. Let us consider
the following map defined on a subset $\Cal V_p$ of the proper domain
$K_p$:
\begin{equation} e_p: \Cal V_p \ni u \mapsto
q=e^{u- K_{p}(u)}\cdot p \in \emmex, \label{eq:ep} \end{equation}
where $ K _{p}(u) =\log \Emu p {e^u} = \log G_p(u)$ is the cumulant
generating functional computed at $u$.

This mapping is one--to--one because $u$ is centered. According to
\formula{chart} and
\formula{patch} we shall denote by $\Cal U_p$ the image of the
mapping and by $s_{p}$ its inverse on $\Cal U_p$.
Such an inverse, $s_p:\Cal U_p\to \Cal V_p$ is easily computed, for
$q\in {\Cal U}_{p}$, as
\begin{equation} \label{eq:expcor}
 s_{p}(q)= \log {q\over p} -
\Emu{p}{\log{q\over p}}. \label{eq:sp}
\end{equation}
{\em The functions $s_p$, $p \in \emmex$, will be the coordinate
mappings
of our manifold in the sense that, locally around each $p \in
\emmex$, each $q \in \Cal U_p$ will be ``parameterized'' by its
{\em centered log--likelihood}}.

Let us compute now the change-of-coordinates formula: if $p_1$ and
$p_2$ are two points in $\emmex$ such that ${\Cal U}_{p_1} \cap
{\Cal U}_{p_2}\neq\emptyset $, then the composite (transition)
mapping

$$ s_{p_2}\circ e_{p_1}: s_{p_1}({\Cal U}_{p_1}\cap {\Cal U}_{p_2})
\rightarrow s_{p_2}({\Cal U}_{p_1}\cap {\Cal U}_{p_2}) $$ simplifies
to $$ s_{p_2}\circ e_{p_1}(u) =u+\log {p_1\over p_2} -
\Emu{p_2}{u+\log {p_1\over p_2}} $$ where the algebraic
computations are done in the space of $\mu$--classes of measurable
functions and the expectation is well defined as long as $\Cal
U_{p_1}\cap \Cal U_{p_2}\ne\emptyset$ because this implies
$u+\log\frac{p_1}{p_2}\in \Cal V_{p_2}$.

\begin{theorem} The collection of pairs $\{(\Cal U_p, s_p): p\in
\emmex\}$ is an affine $C^{\infty }$--atlas on $\emmex$. The induced
topology on sequences is equivalent to $e$--convergence.
\label{th:expoman}\end{theorem}

\begin{definition}[Exponential manifold] The exponential manifold is
the manifold defined by the property in theorem \ref{th:expoman} on
the set $\emmex$. \end{definition}

The manifold structure we have defined is a special one: many other
types of atlases have been suggested in the literature, in
particular the mixture coordinates and the so-called Amari's
imbeddings described in \cite{AMR82}.

%%, see also section \ref{projfm}.
In the infinite dimensional case those different geometric structures
are not equivalent to the exponential manifold, but in some
restricted sense they are, because they induce the same manifold
structure on finite dimensional sub-manifolds (i.e. parametric
statistical manifolds).
%In fact, in section \ref{projfm} we show that the manifold structure
%induces exactly the Fisher information metric on a parametric
%submanifold of densities.
%Such geometric structures are also used to introduce different
%connections, see again the works by Amari and the book by Murray and
%Rice \cite{MRR93}.

The maximal exponential model already defined  has a
precise place in the general framework. In fact the
maximal exponential model $\Cal E_p$ is the
connected component containing $p$ of the exponential manifold
$\emmex$.

In the previous works on the differential geometric approach to
nonlinear
filtering and to the finite dimensional approximation of the
Fokker--Planck equation (Brigo, Hanzon and LeGland \cite{BrHaLe},
\cite{BrHaLe2} and Brigo
\cite{Brigo4}) we used the $L_2$ structure to project the
Kushner--Stratonovich or the Fokker--Planck equation onto a
finite dimensional exponential manifold of densities.
This procedure uses the map $p \mapsto \sqrt{p}$ from positive
densities to their square roots as a tool which allows the
$L_2$ structure to enter the picture. Although this is useful to
perform computations, and even if this approach yields the same
finite dimensional approximation as in the case where one projects
according to the exponential manifold structure discussed here
(compare formulae given in section \ref{projfm} with formulae obtained
via the $L_2$ structure given in \cite{Brigo4}),
we notice that this map cannot be used to define a manifold structure.
It does not yield charts. This is due to the fact that any open set
of $L_2$ contains functions which are negative in a set with positive
measure. Then we see that a chart should map open sets in the manifold
onto open sets in $L_2$, but these open sets would contain the
functions described above, and hence they could not be contained in
any set of square roots of densities (which are positive everywhere).
This is why the space of square roots of densities cannot have
a manifold structure based on $L_2$. In \cite{BrHaLe2} this problem is
bypassed by defining a parametric exponential enveloping manifold.
Here we use the exponential manifold structure to render the procedure
rigorous in an infinite dimensional context.

Now we show the properties of the map from ${\cal M}$ to $L_2$
defined by $R : \ p \mapsto \sqrt{p}$.

\begin{proposition} The mapping
$$R : \Cal M \ni p \mapsto \sqrt p \in L^2(\mu)$$
is $C^\infty$. If the tangent space is identified with $B_p$
then its tangent map is
$$T_p R (v) = \frac 12 R(p) v.$$
In particular the tangent map is surjective at any $p$.
\end{proposition}
\proof Let us fix a density $p_0$ and consider
the coordinate form of the map $R$;
it is defined from ${\cal V}_{p_0}$ to $L_2(p_0)$ by
$H_{p_0}(u) := \sqrt{e_{p_0}(u)}$. By direct computation one
obtains the form of the directional derivative
\begin{displaymath}
{d \over dv} H_{p_0}(u) =  H_{p_0} (u) (\half [v - D K_{p_0}(u)v])
   =  H_{p_0}(u) \half [v - E_{e_{p_0}(u)} v],
\end{displaymath}
and the norm of the differential operator
\begin{displaymath}
\Vert D \ H_{p_0}(u) \ v \Vert_2^2 =
   \inverse{4} E_{e_{p_0}(u)} \{( v - E_{e_{p_0}(u)} v )^2 \}
   = \inverse{4} D^2 K_{p_0}(u)(v,v),
\end{displaymath}
which shows that $H_{p_0}$ is differentiable at $u$ since $K_{p_0}$
is infinitely  Fr\'echet differentiable. The coordinate-free form of the
differential is $T_p R v = \half \sqrt p v$ if we identify the
tangent space at $p$ with $B_p$.
Note that the $L_2$ norm of this first derivative represents a
variance, so that the derivative operator is one-to-one.

In a similar way the other derivatives of $H_{p_0}$ can be computed as:
\begin{displaymath}
 D^n \ H_{p_0}(u) (v_1, \dots, v_n) =
2^{-n} H_{p_0} (u) (v_1 - E_{e_{p_0}(u)} v_1, \dots ,
v_n - E_{e_{p_0}(u)} v_n).
\end{displaymath}
If one computes the $L_2$ norm of this derivative one finds
easily that the norm is bounded and the differentiability of any
order follows.
\QED
Hence the mapping $p \mapsto \sqrt p$ is what we
call a regular parametrization. From the exponential
coordinates $u$ we deduce
$H_{p_0}(u)$ which can be differentiated. Yet $H_{p_0}$, although
constituting a parametrization, does not define coordinates.
This is due to the fact that we are working in infinite dimension and
a regular parametrization is not necessarily a chart.

A basic object of the manifolds theory is the tangent bundle. In the
case of the exponential manifold it has been remarked from the very
beginning of the geometrical theory that there is a very natural
identification between the tangent vectors and the exponential
one-dimensional models around a point $p$. In fact each
differentiable curve in $\emmex$, i.e. each one-dimensional
statistical model  $p(t)$, $t\in I\subseteq R$, such that $p(0) = p$,
has a tangent model of the exponential form $e^{tu - K_p(tu)} \cdot
p$. This can be rephrased by saying that
these exponential one-dimensional models seem to play the role of
straight lines.

\begin{definition}[Tangent space] The tangent space $T_p$ at $p$ of
the exponential manifold on $\emmex$ is the set (indexed by $u$) of the
one-dimensional exponential models

\begin{displaymath} e^{tu - K_p(tu)} \cdot p,\quad t\in\reals,\quad
u\in B_p \quad . \end{displaymath}

Usually we will identify the tangent exponential model $e^{tu -
K_p(tu)} \cdot p$ with its {\em score}

$$\left. \frac d {dt} \left(tu -
K_p(tu)\right) \right|_{t=0} = u \in B_p \quad .$$

The tangent space
inherits the structure of Banach space from $B_p$. \end{definition}
\begin{definition}[Sub-manifold, sub-model\label{def:submani}]
Let $\Cal N$ be a subset of the exponential manifold $\emmex$ and,
for each density $p \in {\Cal N}$, let $V_p^1$ and $V_p^2$ be closed
subsets of $B_p$, such that there exist:
\begin{enumerate}
\item
 a linear invertible and bi-continuous
mapping between  $B_p$ and some direct sum $V_p^1 + V_p^2$.
That is  $V_p^1$ and $V_p^2$ {\em split} in $B_p$.
\item
a chart on a neighbourhood $\Cal W_p$ of $p$:
\begin{displaymath}
\sigma_p: \Cal W_p \to V_p^1 + V_p^2,
\end{displaymath}
where  $\sigma_p$
maps $\Cal W_p$ onto the product of to open sets $\Cal V_p^1 \times \Cal
V_p^2$ and ${\Cal N} \cap \Cal W_p$ onto $\Cal V_p^1 \times \set 0$.
\end{enumerate}
We
will say that ${\Cal N}$ is a {\em sub-model} or a {\em
sub-manifold} of the exponential manifold $\emmex$. \end{definition}

A sub-manifold is a manifold defined by the restricted maps.
For a list of examples see \cite{PSR94a}. Our basic example is a
finite dimensional exponential family
\begin{eqnarray*}
  && EM(c) = \{p(\cdot,\theta), \theta \in \Theta\},\\
  && p(\cdot,\theta) := \exp[\theta^T c(\cdot) - \psi(\theta)],
\end{eqnarray*}
where $c=(c_1, \ldots, c_n)$, and $\Theta$ is a convex open set
in $\realvectors n$.
In this case the local representation at $p(\cdot, \theta_0)$ is
\begin{displaymath}
p(\cdot,\theta) =
\exp[(\theta - \theta_0)^T [c(\cdot) - \psi'(\theta_0)] -
(\psi(\theta) - \psi(\theta_0)) +
(\theta - \theta_0)^T \psi'(\theta_0)] p(\cdot, \theta_0),
\end{displaymath}
and the relevant splitting is
\begin{eqnarray*}
&& V^1_{p(\theta_0)} = \span
\set{c_i - {\partial \over \partial \theta_i} \psi(\theta_0)}\\
&& V^2_{p(\theta_0)} = \set{u \in B_{p(\theta_0)} :
\Emu {p(\theta_0)}{uc_i} =0, i=1, \ldots, n} .
\end{eqnarray*}
\section{Evolution of marginal laws of a diffusion process}
On the complete probability space $(\Omega,{\cal F},P)$ let us
consider a stochastic process $\{X_t, t \ge 0\}$ of diffusion type.
%adapted to a filtration  $\{{\cal F}_{t}, t \ge 0 \}$.
Let the dynamic equation describing $X$ be
of the following form
\begin{eqnarray*}
   dX_t =  f_t(X_t) dt + \sigma_t(X_t) d W_t,
\end{eqnarray*}
where $\{W_t, t\ge 0\}$ is a standard Brownian motion independent
of the initial condition $X_0$.
The equation above is an It{\^o} stochastic differential
equation.
In the following derivation,
we treat the scalar case.
The following set of assumptions will be in force throughout
the paper.

\begin{itemize}
\item[(A)] Initial condition:~ We assume that the initial
state $X_0$ has a density $p_0$
w.r.t.\ the Lebesgue measure on $\reals$,
%and has finite moments of any order,
with $p_0$ almost surely positive.
%
%and we make the following assumptions
%on the coefficients $f_t$, $a_t := \sigma_t^2\,$:
%\begin{itemize}
   \item[(B)] Local strong existence:~$f \in C^{1,0}$, $a \in  C^{2,0}$,
   which means that $f$ is once continuously differentiable
   wrt $x$ and continuous wrt  $t$
   and $a$ is twice continuously differentiable
   wrt $x$ and continuous wrt $t$.
   This assumptions imply in
   particular local Lipschitz continuity.

   \item[(C)] Non--explosion~:
there exists $K > 0$ such that
\begin{eqnarray*}
  2 x f_t(x) + a_t(x) \leq K\, (1+\vert x \vert^2),
\end{eqnarray*}
for all $t\geq 0$, and for all $x\in \reals$.
\end{itemize}

Under assumptions~(A),~(B) and~(C) there exists
a unique solution $\{X_t\,,\,t\geq 0\}$ to the state equation,
see \cite{StroVarh}, theorem 10.2.1 with
$\phi(x,t) = x^2$.

%Under additional assumptions on the coefficients
%(see e.g. \cite{FrieA} or \cite{StroVarh})
\begin{itemize}
\item[(D)]
We assume that the law of $X_t$ is absolutely continuous and
its density
$p_t(x)$ has regularity $C^{2,1}$ and satisfies
the Fokker--Planck equation (FPE):
\begin{eqnarray} \label{eq:FP}
\frac{\partial p_t}{\partial t} = {\cal L}_t^\ast p_t,
\end{eqnarray}
where the backward diffusion operator ${\cal L}_t$
is defined by
\begin{displaymath}
   {\cal L}_t = f_t\,
   \frac{\partial}{\partial x} + \half
   a_t \frac{\partial^2}{\partial x^2},
\end{displaymath}
and its dual (forward) operator is given by
\begin{displaymath}
   {\cal L}^\ast_t p = -
   \frac{\partial}{\partial x} (f_t p) + \half
   \frac{\partial^2}{\partial x^2} (a_t p).
\end{displaymath}
We assume also $p_t(x)$ to be positive for all $t \ge 0$ and
almost all $x \in \reals$.
\end{itemize}
Assumption (D) holds for example under conditions given by
boundedness of the coefficients $f$ and $a$ plus
uniform ellipticity of $a_t$, see \cite{StroVarh} theorem 9.1.9.
Different conditions are also given in \cite{FrieA},
theorem 6.4.7.
Now we rewrite equation (\ref{eq:FP}) in the exponential
coordinates (\ref{eq:expcor}).
%
%We assume that the curve $t \mapsto p_t$ in ${\cal M}$
%is differentiable. In such a case we compute from equation
%(\ref{eq:FP}) the representation of the tangent vectors
%in exponential coordinates.
Consider as local reference density the
solution $p_t$ of FPE at time $t$. We are now working around $p_t$.
Consider a curve around $p_t$ corresponding to the
solution of FPE around time $t$ expressed in $B_{p_t}$
coordinates:
\begin{eqnarray*}
(-\epsilon,\epsilon) &\rightarrow& {\cal V}_{p_t} \\
h &\mapsto& s_{p_t}(p_{t+h})=: u_h.
\end{eqnarray*}
The function $u_h$ represents the expression in coordinates of
the density
\begin{eqnarray} \label{coordin}
p_{t+h} = \exp[u_h - K_{p_t}(u_h)] p_t =: e_h p_t.
\end{eqnarray}
Now consider FPE around $t$, i.e.
\begin{eqnarray*}
\frac{\partial p_{t+h}}{\partial h} = {\cal L}_{t+h}^\ast p_{t+h}.
\end{eqnarray*}
Substitute (\ref{coordin}) in this last equation in order to obtain
\begin{eqnarray*}
\frac{\partial e_h p_t}{\partial h} = {\cal L}_{t+h}^\ast (e_h p_t).
\end{eqnarray*}
Write
\begin{eqnarray*}
\frac{\partial e_h}{\partial h} = \frac{{\cal L}_{t+h}^\ast (e_h
p_t)}{p_t}
\end{eqnarray*}
and set $h=0$, since we are concerned with the behaviour in $t$.
Notice that $e_0=\exp[u_0-K_{p_t}(u_0)]=\exp(0)=1$,
and that
\begin{eqnarray*}
\frac{\partial e_h}{\partial h}|_{h=0} = \{e_h
\frac{\partial [u_h - K_{p_t}(u_h)]}{\partial h}\}|_{h=0}
=\frac{\partial [u_h - K_{p_t}(u_h)]}{\partial h}|_{h=0}.
\end{eqnarray*}
Moreover, by straightforward computations (write explicitly
the map $K_{p_t}$, use
$u_h=s_{p_t}(p_{t+h})$ and differentiate wrt $h$ under the expectation
$E_{p_t}$) one verifies
\begin{eqnarray*}
\left. \frac{\partial K_{p_t}(u_h)}{\partial h} \right |_{h=0}=0,
\end{eqnarray*}
so that
\begin{equation} \label{equation-u}
\left. \frac{\partial u_h}{\partial h} \right |_{h=0}
= \frac{{\cal L}_t^\ast p_t}{p_t}
\end{equation}
is the formal representation in exponential coordinates
of the tangent vector in $p_t$.
Notice that, again by straightforward computations,
\begin{eqnarray} \label {alpha-def}
   \alpha_t := \alpha_t(p_t) = \frac{{\cal L}_t^\ast p_t}
          {p_t}
   &=&  - f_t\, \frac{\partial}{\partial x}(\log p_t)
    - \frac{\partial f_t}{\partial x} +
   \\ \nonumber \\ \nonumber
   &+& \half  [\,
   a_t\,
   \frac{\partial^2}{\partial x^2}(\log p_t)
   + a_t\, (\frac{\partial}{\partial x}(\log p_t))^2\,+
   \\ \nonumber \\  \nonumber
   &+& 2\, \frac{\partial a_t}{\partial x}\,
   \frac{\partial}{\partial x}(\log p_t)
   + \frac{\partial^2 a_t}{\partial x^2} \,]\ .
\end{eqnarray}
Summarizing: consider the curve expressing FPE around $p_t$ in
$B_{p_t}$ coordinates. Its tangent vector is given
by $\alpha_t$.
Under suitable assumptions on the coefficients $f_t$ and $a_t$
the function
$\alpha_t$ belongs to $B_{p_t}$, according to the
convention that locally identifies the tangent space of normed spaces
with the normed space itself.
To render the computation not only formal
we need $\alpha_t$ to be really a tangent vector for our
manifold structure. This in turn requires the curve
$t \mapsto p_t$ to be differentiable.
Below we give a regularity result
expressing a condition under which this happens and whose proof
is immediate.
Moreover, we give a condition which can be used to check
whether the evolution stays in a given submanifold.
\begin{proposition}
[Regularity and finite dimensionality
of the solution of FPE] \label{th:reg}
\hspace{2cm}
\begin{itemize}
\item[(i)] If the map $t \mapsto p_t$ is differentiable
in the manifold ${\cal M}$ then $\alpha_t$ given in
eq. (\ref{alpha-def}) is a tangent vector.
\item[(ii)] If the map
$t \mapsto \alpha_t$ is continuous
at $t_0$ into $L^{\cosh \cdot -1}$,
then $t \mapsto p_t$ is differentiable at $t_0$
as a map into ${\cal M}$.
\item[(iii)]
Let be given a submanifold ${\cal N}$ such
that $p_0 \in {\cal N}$. If the previous condition
is satisfied and
\begin{displaymath}
\frac{{\cal L}_t^\ast p}{p}
\end{displaymath}
is tangent to ${\cal N}$ at $p$ for all $p \in {\cal N}$,
then $p_t$ evolves in ${\cal N}$.
\end{itemize}
\end{proposition}

Sufficient conditions under which condition (ii) in the
proposition happens to be true
are given by boundedness for all possible $T>0$ of $f$, $\partial_x f$,
$a$, $\partial_x a$, $\partial^2_{xx} a$ in $[0 \ \ T] \times I\!\!R$
plus classical assumptions given in Stroock and Varadhan
\cite{StroVarh} theorem 9.1.9,
or Friedman \cite{FrieA}, theorem 6.4.7, ensuring existence of a
regular solution or Fokker--Planck equation (as required in (D)).
This follows from the fact that if $\alpha_t(x)$ is continuous
and bounded in both $t$ and $x$, then it is continuous as
a map $t \mapsto \alpha_t$ from $[0 \ \ T]$ to $L^{\cosh \cdot -1}$.

In the following we give examples where this proposition
applies. Some of them are obtained from \cite{Brigo4}
where the detailed derivation is given.
\begin{example}[Linear case]
If $f_t(x) = F_t x$ for all $t \ge 0, \ x \in \reals$
($f$ linear in $x$) and if $a_t(x) = A_t$ for all
 $t \ge 0,\ x \in \reals$ ($a$ does not depend on $x$) and if
 finally $p_0 \sim {\cal N}(m_0,Q_0)$ then it is known that
$p_t \sim {\cal N}(m_t,Q_t)$
where $m_t = m_0 \exp{\int_0^t F_s ds} $ and
$Q_t$ is the (unique) positive solution of the
(scalar Lyapunov) equation
\begin{eqnarray*}
\dot{Q_t} = 2 F_t Q_t + A_t,
\end{eqnarray*}
with initial condition $Q_0$ given.
%
%{\cal N}(\mu_0 \exp(F\ t),
%[(2F \sigma_0^2 + A) \exp(2Ft)-A] /2F)$.
Consider now a generic Gaussian density $p \sim {\cal N}(m,Q)$
and compute
\begin{eqnarray} \label{eq:exlin}
(\frac{{\cal L}_t^\ast p}{p})(x) =
  (\frac{F_t}{Q} + \frac{A_t}{2 Q^2}) x^2 -
  (\frac{F_t m}{Q} + \frac{A_t m}{Q^2}) x +
  \frac{A_t m^2 }{2 Q^2} - F_t - \frac{A_t}{2 Q}.
\end{eqnarray}
When applied to $p_t$, the previous formula yields
$\alpha_t$:
\begin{eqnarray*}
\alpha_t
 = (\frac{F_t}{Q_t} + \frac{A_t}{2 Q_t^2}) x^2 -
  (\frac{F_t m_t}{Q_t} + \frac{A_t m_t}{Q_t^2}) x +
  \frac{A_t m_t^2 }{2 Q_t^2} - F_t - \frac{A_t}{2 Q_t},
%
%=\frac{{\cal L}^\ast p_t}{p_t} = -\partial_t \sigma_t
%    - \frac{(x-\mu_t)}{\sigma_t^3} \
%      [ - \sigma_t \dot{\mu_t} - (x - \mu_t) \dot{\sigma_t} ].
\end{eqnarray*}
where $m_t$ and $Q_t$ have been defined above.

In this case the previous proposition applies.
First, one sees that $t \mapsto \alpha_t$ is indeed
continuous at any $t_0$ in $L^{\cosh \cdot - 1}$.
Secondly, one can deduce already from (\ref{eq:exlin})
{\em without solving the Fokker--Planck equation} that the solution
will have a Gaussian density. Indeed, one can easily check
that the tangent space to the Gaussian submanifold
of ${\cal M}$ expressed in $B$ coordinates
contains the function space $\span\{1,x,x^2\}$.
Since by expression (\ref{eq:exlin}) we see that
$({\cal L}_t^\ast p)/p$ lies in $\span\{1,x,x^2\}$
for all $p$ in the Gaussian submanifold, we deduce
that the solution of the Fokker--Planck equation will evolve in
the Gaussian submanifold.
\end{example}
\begin{example} [Nonlinear diffusions with
unit variance Gaussian law]
Let be given a diffusion coefficient $\sigma_t(x)$ satisfying
assumptions (B) and assumption (C) when the
drift vanishes, i.e. when $f=0$
(we set as usual $a := \sigma^2$).
In \cite{Brigo4} it is shown that by defining the drift
\begin{eqnarray*}
&& f_t(x) := \half \frac{\partial a_t}{\partial x}(x) +
          \half a_t(x)[ kt - x ] + k,
\end{eqnarray*}
the Fokker--Planck equation for the density of the solution of
the stochastic differential equation
\begin{eqnarray*}
   d X_t = f_t(X_t) dt + \sigma_t(X_t) dW_t \ \ \
   X_0 \sim {\cal N}(0,1),
\end{eqnarray*}
is solved by $p_t \sim {\cal N}(kt,1)$ for all possible
diffusion coefficients $\sigma_t(x)$.
Here the solution of the Fokker--Planck equation evolves in a
submanifold of ${\cal M}$ given by Gaussian densities
with unit variance. Actually, the mean of $p_t$
evolves linearly in time and the variance is fixed to one.
Note that in this case
\begin{eqnarray*}
\alpha_t = \partial_t \ \log p_t = k (x - kt),
\end{eqnarray*}
and the curve $t \mapsto \alpha_t$ is clearly continuous
at any $t_0$ in $L^{\cosh \cdot -1}$. One might check a priori that
if a given $p$ belongs to the submanifold of Gaussian densities with
unit
variance, then $({\cal L}_t^\ast p)/p$ belongs to the tangent space
of this submanifold if the mean is given by $kt$.
Indeed, we are considering the family
\begin{eqnarray*}
  p(x,\theta) = \frac{1}{\sqrt{2 \pi}} \exp[-\half(x - \theta)^2] \sim
{\cal N}(
\theta,1),
  \ \ \theta \in \reals,
\end{eqnarray*}
and its tangent space expressed in $B_{p(\cdot,\theta)}$ coordinates,
$\span\{x-\theta\}$.
Let us compute
\begin{eqnarray*}
\alpha_{t,\theta}(x) =  \half (\partial_x a_t(x)) \ ( kt - \theta ) + \
    \half a_t(x) (x - \theta) (kt - \theta)
    \ + \ k (x - \theta).
\end{eqnarray*}
Under reasonable assumptions on $a$, this function
belongs to the tangent space $\span\{x-\theta\}$ if and only if
$\theta = kt$. We have been able to check that the density
of the diffusion $X$ evolves according to $p_t \sim {\cal N}(kt,1)$
without solving the Fokker--Planck equation.
\end{example}
>From examples given in \cite{Brigo4} one can construct
  other nonlinear cases where the above proposition
applies.
\section{Projection of the Fokker--Planck equation}
\label{projfm}
In reaching equation (\ref{equation-u}) we assumed
implicitly a few facts.
We are assuming that there always exists a neighborhood of
$h=0$ such that in this neighborhood $p_{t+h} \in {\cal U}_{p_t}$.
Conditions under which this happens will be examined in the future.
We only remark that when projecting on a finite dimensional
exponential manifold, these conditions are not necessary for the
projected equation to exist and make sense, see below.
Neither we need equation (\ref{equation-u}) to have a solution
to obtain existence of the solutions of the projected equation.
Now we shall project this equation on a finite dimensional
parametrized exponential manifold.
First notice that according to proposition \ref{th:reg}
the law $p_t$ remains finite dimensional if
there exists a finite dimensional parametrized submanifold of
${\cal M}$,
\begin{eqnarray*}
   S = \{p(\cdot,\theta), \theta \in \Theta\},
\end{eqnarray*}
such that the corresponding tangent vectors
$\alpha_t(p(\cdot,\theta))$ of the FPE are in the tangent space
of this finite dimensional submanifold.
We take the set $\Theta$ open in $\reals^m$.
If we look for a finite dimensional {\em exponential} family,
we can select a submanifold
\begin{eqnarray*}
  && EM(c) = \{p(\cdot,\theta), \theta \in \Theta\},\\ \\
  && p(\cdot,\theta) := \exp[\theta^T c(\cdot) - \psi(\theta)].
\end{eqnarray*}
We will assume the following on the family $EM(c)$ (see \cite{BrHaLe},
\cite{BrHaLe2} for other more specific assumptions):
\begin{itemize}
\item[(E)] We assume $c \in C^2$.
%and $c(\cdot)$ to have at most
%polynomial growth.
\end{itemize}
Notice that tangent vectors around a point $p(\cdot,\theta_t)$
of a generic curve $h \mapsto p(\cdot,\theta_{t+h})$ on $EM(c)$,
are now obtained according to (after straightforward computations
and the chain rule)
\begin{eqnarray} \label{tang-vect}
\frac{\partial s_{p(\cdot,\theta_t)}(p(\cdot,\theta_{t+h})) }
{\partial h}|_{h=0} =
\sum_{i=1}^{m}[c_i(\cdot) - E_{\theta_t}c_i] \dot{\theta}_t^i.
\end{eqnarray}
As a consequence, the tangent space at $\theta$ is given by
\begin{displaymath}
T_\theta EM(c) =
\span\{c_1(\cdot)-E_{\theta}c_1,\ldots,c_m(\cdot)-E_{\theta}c_m\},
\end{displaymath}
where $E_{\theta}\{\phi\}:= \int \phi(x) p(x,\theta) dx$.
Consider the following inner product in
$T_\theta EM(c)$ (and more in general in $B_{p(\cdot,\theta)}$):
\begin{eqnarray*}
\langle v_1,v_2 \rangle_\theta := E_\theta[v_1\ v_2],
\hspace{1cm}
v_1,v_2 \in B_{p(\cdot,\theta)}.
\end{eqnarray*}
Consider the quantities
\begin{displaymath}
g(\theta)_{ij} :=
  \langle
  c_i(\cdot)-E_{\theta}c_i,\
  c_j(\cdot)-E_{\theta}c_j
  \rangle_\theta, \ \ i,j=1,\ldots,m.
\end{displaymath}
Notice that the matrix $g(\theta)$, expressing the inner products of
tangent vectors in $T_\theta EM(c)$, is nothing else than
the traditional Fisher information matrix
\begin{eqnarray*}
(\ E_{\theta}\ [\partial_{\theta_i}  \log p(\cdot,\theta) \ \
    \partial_{\theta_j} \log p(\cdot,\theta)] \ )_{i,j=1,\ldots,m}
\end{eqnarray*}
for the family $EM(c)$ (see also \cite{BrHaLe}, \cite{BrHaLe2}).
Now define for all $\theta\in\Theta$
the orthogonal projection
\begin{eqnarray*}
  && \Pi_\theta : B_{p(\cdot,\theta)} \longrightarrow T_\theta EM(c) \\
  && \Pi_\theta [ v ] :=
   \sum_{i=1}^m [ \sum_{j=1}^m g^{ij}(\theta)\;
   \langle v, c_j(\cdot) - E_{\theta}c_j
 \rangle_\theta ]\; ( c_i(\cdot) - E_{\theta}c_i) .
\end{eqnarray*}
A rapid computation involving
duality between ${\cal L}$ and ${\cal L}^\ast$ and standard results
on the normalization constant $\psi(\theta)$ of exponential families
(such as $\partial_{\theta_i} \psi(\theta) = E_{\theta} c_i$)
yields
\begin{eqnarray*}
 {\cal P}_{t,\theta}
  &:=& \Pi_\theta
  [\frac{{\cal L}_t^\ast p(\cdot,\theta)}
          {p(\cdot,\theta)}] =
    E_{\theta}[{\cal L}_t c]^T \  g^{-1}(\theta)
        \ [c(\cdot) - E_{\theta} c],
\end{eqnarray*}
where integrals of vector functions are meant to be applied to their
components.
Note that this map is regular in $\theta$ under reasonable assumptions
on $f,a$ and $c$.
At this point we project equation (\ref{equation-u})
via this projection. By remembering expression (\ref{tang-vect})
for tangent vectors and the above formula for the projection
we obtain the  following ($m$--dimensional)
ordinary differential equation (in vector form)
on the manifold $EM(c)$:
\begin{eqnarray}  \label{PFPE}
[c(\cdot) - E_{\theta_t}c]^T \dot{\theta}_t =
E_{\theta_t}[{\cal L}_t c]^T \  g^{-1}(\theta_t)
\ [c(\cdot) - E_{\theta_t}c].
\end{eqnarray}
It follows immediately the following
ordinary differential equation for the parameters:
\begin{eqnarray} \label{PFPEPAR}
   \dot{\theta}_t
   = g^{-1}(\theta_t)\; E_{\theta_t}\{{\cal L}_t\; c\}.
\end{eqnarray}
Notice that, as anticipated above,
equation (\ref{PFPEPAR})
is well defined and admits locally a unique solution
if  the following condition (ensuring existence of the
norm of $\alpha_t(p(\cdot,\theta_t))$ associated to the inner product
$\langle \cdot\ ,\ \cdot \rangle_{\theta_t}$)
holds:
\begin{eqnarray} \label {alpha-cond}
(F) \hspace{1cm}   &&E_{\theta}\{\alpha_{t,\theta}^2\}<\infty \;\;
                  \forall \theta \in \Theta,
   \\ \nonumber \\ \nonumber
   &&\alpha_{t,\theta} := \frac{{\cal L}_t^\ast p(\cdot,\theta)}
          {p(\cdot,\theta)} =
     - f_t\, \frac{\partial}{\partial x}(\theta^T c)
    - \frac{\partial f_t}{\partial x} +
       \\ \nonumber \\   \nonumber
   &&\hspace{1cm}  + \half  [\,
   a_t\,
   \frac{\partial^2}{\partial x^2}(\theta^T c)
   + a_t\, (\frac{\partial}{\partial x}(\theta^T c))^2\,+
   \\ \nonumber \\  \nonumber
   && \hspace{1cm} + 2\, \frac{\partial a_t}{\partial x}\,
   \frac{\partial}{\partial x}(\theta^T c)
   + \frac{\partial^2 a_t}{\partial x^2} \,]\ .
\end{eqnarray}
We will assume such condition to hold in the following.
Notice that this is a condition on the coefficients $f,a,c$.
We have thus proven the following
\begin{proposition}
[ Projected evolution of the density of an It\^o diffusion]
\label{proj-law}
Assume assumptions (A), (B),(C), (E) and (F) on the coefficients
$f, a$, on the initial condition $X_0$ of the It\^o diffusion $X$,
and on the sufficient statistics
$c_1,\dots,c_n$ of the exponential family $EM(c)$ are satisfied.
Then the projection of Fokker--Planck equation describing the
evolution of $p_t = p_{X_t}$ onto $EM(c)$  reads, in
$B_{p_t}$ coordinates:
\begin{eqnarray*}
[c(\cdot) - E_{\theta_t}c]^T \dot{\theta}_t =
E_{\theta_t}[{\cal L}c]^T \  g^{-1}(\theta_t)
\ [c(\cdot) - E_{\theta_t}c],
\end{eqnarray*}
and the differential equation describing the evolution of the
parameters for the projected density--evolution is
\begin{eqnarray*}
   \dot{\theta}_t
   = g^{-1}(\theta_t)\; E_{\theta_t}\{{\cal L}_t\; c\}.
\end{eqnarray*}
\end{proposition}
Notice that the projected equations exist under conditions
which are more general than conditions for existence of the
solution of the original Fokker--Planck equation.
For more details see Brigo \cite{Brigo4}.

\section{Interpretation of the projected density as density
of a different diffusion}
In this section we shortly expose a problem which was treated in
\cite{Brigo4}.
Consider the projected density $p(\cdot,\theta_t)$,
expressing the projection of the
density--evolution of the one dimensional diffusion $X$ onto the
exponential manifold $EM(c)$.
The question is: Can we define a diffusion
$Y_t$ whose density is the
projected density $p(\cdot,\theta_t)$ ?
If the answer is yes, $Y_t$ is a diffusion
whose density evolves in a finite
dimensional exponential manifold assigned a priori (for example
Gaussian).
In order to proceed,
define a diffusion
\begin{equation}
dY_t = u_t(Y_t)dt + \sigma_t(Y_t) dW_t, \ \ Y_0 = X_0,
\end{equation}
with the same diffusion coefficient as $X_t$. We shall try to define
the drift $u$ in such a way that the density--evolution of $Y_t$
coincides
with $p(\cdot,\theta_t)$.
Call ${\cal T}_t$ the backward differential operator of $Y_t$:
\begin{displaymath}
   {\cal T}_t = u_t \,
   \frac{\partial}{\partial x} + \half
   a_t \frac{\partial^2}{\partial x^2}.
\end{displaymath}
Consider the right hand sides of (\ref{equation-u}) and (\ref{PFPE}).
Clearly, the density of $Y_t$ coincides with $p(\cdot,\theta_t)$
if
\begin{displaymath}
\frac{{\cal T}^\ast p(\cdot,\theta_t)}
   {p(\cdot,\theta_t)} =
E_{\theta_t}[{\cal L}c]^T \  g^{-1}(\theta_t)
\ [c(\cdot) - E_{\theta_t}c]
\end{displaymath}
which we can rewrite as
\begin{displaymath}
{\cal T}^\ast p(\cdot,\theta_t) = {\cal P}_{t,\theta_t} \
p(\cdot,\theta_t).
\end{displaymath}
By simple calculations one can rewrite the above equation
as the following PDE for $u$, where we do not expand the second
partial derivative of $a_t p(\cdot,\theta)$:
\begin{eqnarray*}
  \frac{\partial u_t}{\partial x} +
  \theta_t^T \frac{\partial c}{\partial x} u_t =
   \frac{1}{2\, p(\cdot,\theta_t)}
   \frac{\partial^2}{\partial x^2}(a_t p(\cdot,\theta_t)) -
   {\cal P}_{t,\theta_t}
\end{eqnarray*}
Call ${\cal B}_{t,\theta_t}$ the right hand side of such equation.
A solution is given by
\begin{displaymath}
u^\ast_t(x):= \exp[-\theta_t^T c(x)] \int_{-\infty}^x
               {\cal B}_{t,\theta_t}(y) \exp[\theta_t^T c(y)] dy,
\end{displaymath}
as one can verify immediately by substitution.
Straightforward calculations yield
\begin{eqnarray} \label{ustar}
 u^\ast_t(x) &:=& \frac{1}{p(x,\theta_t)}
  \int_{-\infty}^{x} [
    \frac{\partial^2_{xx} (a_t(y)p(y,\theta_t))}{p(y,\theta_t)}
    - \Pi_{\theta_t}
      \{\frac{\partial^2_{xx} (a_t(y)p(y,\theta_t))}{p(y,\theta_t)} \} +
      \\ \nonumber \\ \nonumber
      && \hspace{2cm}
        +\Pi_{\theta_t} \{\frac{\partial_x (f_t(y) p(y,\theta_t))}
                {p(y,\theta_t)}\} ]
      p(y,\theta_t)  dy = \\ \nonumber \\ \nonumber
  &=&\half \frac{\partial a_t}{\partial x}(x) +
          \half a_t(x) \theta_t^T  \frac{\partial c}{\partial x}(x)
          + \\ \nonumber \\ \nonumber
&& - E_{\theta_t}\{{\cal L}_t c\}^T g^{-1}(\theta_t)
                 \int_{-\infty}^x (c(y) - E_{\theta_t}c)
                 \ \exp[\theta_t^T (c(y) - c(x))] dy.
\end{eqnarray}
>From this last equation one sees that under condition
(\ref{alpha-cond}) and under the assumption that densities of $EM(c)$
are integrable, the above integral always exists.

We have thus proven the following
\begin{proposition} [Interpretation of the projected density--evolution]
Assume assumptions (A), (B), (C), (E) and (F) on the coefficients
$f, a$ and on the initial condition $X_0$
of the It\^o diffusion $X$ and on the sufficient statistics
$c$ of the exponential family $EM(c)$ are satisfied.
Let $p(\cdot,\theta_t)$ be the projected density evolution, according to
proposition \ref{proj-law}. Define
\begin{eqnarray*}
d Y_t &=& u^\ast_t(Y_t) dt + \sigma_t(Y_t) dW_t, \\ \\
 u^\ast_t(x) &:=& \half \frac{\partial a_t}{\partial x}(x) +
          \half a_t(x) \theta_t^T  \frac{\partial c}{\partial x}(x)
          + \\ \nonumber \\ \nonumber
&& - E_{\theta_t}\{{\cal L}_t c\}^T g^{-1}(\theta_t)
                 \int_{-\infty}^x (c(y) - E_{\theta_t}c)
                 \ \exp[\theta_t^T (c(y) - c(x))] dy.
\end{eqnarray*}
Then $Y$
is an It\^o diffusion whose density--evolution coincides with the
projected density--evolution
$p(\cdot,\theta_t)$ of $X_t$ onto $EM(c)$.
\end{proposition}

\section{Further research on convergence of the projected density
towards the original one}
Now we consider the problem of the convergence
of the
projected  density--evolution to the evolution of the limit diffusion.
The idea is to consider a sequence of
nested finite dimensional families and to check what happens when
the dimension of the family on which the equation is projected tends
to infinity.
The problem we shall investigate in the future is the following.
Suppose we can define a sequence of families
in the following way:
\begin{description}
\item[(G)] we are given a sequence of functions
 $(c_j)_{j \in I\!\!N}$.
Call $c^m := \{c_1 \ c_2 \ \ldots c_m\}$, and assume that
for all $m$ the family $EM(c^m)$ is a finite dimensional
exponential manifold satisfying assumptions (E) and (F).
Call $p(\cdot,\theta^m_t)$ the density coming from projection
of Fokker--Planck equation onto $EM(c^m)$.
\end{description}
As we saw in the preceding section (see (\ref{ustar})) this is
also the density of a diffusion process $Y^m$ with the same initial
condition, the same diffusion coefficient and drift given by
\begin{eqnarray} \label{ustar-m}
 u^m_t(x) &:=& \frac{1}{p(x,\theta^m_t)}
  \int_{-\infty}^{x} [
    \frac{\partial^2_{xx} (a_t(y)p(y,\theta^m_t))}{p(y,\theta^m_t)}
    - \Pi_{\theta^m_t}
      \{\frac{\partial^2_{xx} (a_t(y)p(y,\theta^m_t))}{p(y,\theta^m_t)}
\}
       + \\ \nonumber \\ \nonumber
      && \hspace{2cm}
       +\Pi_{\theta^m_t} \{\frac{\partial_x (f_t(y) p(y,\theta^m_t))}
                {p(y,\theta^m_t)} \} ]
      p(y,\theta^m_t)  dy.
\end{eqnarray}
Hence, if we prove that $Y^m$ converges in law towards the
original diffusion $X$, we have a first convergence result
of the projected evolution towards the original one.
Now, it is well known that since all the diffusions $Y^m$ share
the same diffusion coefficient, under reasonable assumptions on
$\sigma$ the sequence of the laws of $(Y^m)_m$ is relatively
compact in the space of processes with continuous trajectories
(see for example Stroock and Varadhan \cite{StroVarh} or Bafico and
Pistone
\cite{BafiPist}).
Moreover, if the drifts $u^m$ weakly converge to $f$,
the law of $Y^m$ will converge to the law of $X$.
By looking at expression (\ref{ustar-m}), one sees intuitively that
this should happen. Indeed, assume $p(\cdot,\theta^m)$
admits a limit $\bar{p}(\cdot)$ when $m$ tends to infinity,
and suppose the projection tends to be exact when
$m$ tends to infinity. Then in formula  (\ref{ustar-m})
replace $p(\cdot,\theta^m)$ by $\bar{p}$ and eliminate
the projection operators. The expression for $u_t(x)$ simplifies
to $f_t(x)$, so that we have pointwise convergence of the drifts
$u^m$ towards $f$ and we are done.
Of course, one needs to make the above idea precise, to
show reasonable choices of $c^m$ and to prove the result
rigorously. These problems will be investigated in the next future.
%
%Let $p^m$ be a density in $EM(c^m)$ for all $m$.
%Assume $c^m$ are in $B_{p^m}$ for all $m$.
%which is total in $L_0^2{p\cdot \lambda}$; ? ? ? ?
%
%Let be given a density $p$ in
%$\Cal M(\realvectors d, \Cal B(\realvectors d), \lambda)$.
%We will study
%the evolution of the density of the diffusion $(X_t)$ in the
%neighborhood of the initial density $p$. Assume the following:
%
%
%
%$EM(c^m)$ is the (finite-dimensional) exponential manifold defined by
%the mapping $\Theta_n \ni \theta^n \mapsto e_p(\theta_n \cdot c^n)$,
%where $\Theta_n$ is the interior of the proper domain of $\theta^n
%\mapsto G_{p_0}(\theta^n \cdot C^n)$.
%Let $p^m(\cdot,\theta_t^m)$ be the projection
%of the law of $X_t$ onto $EM(c^m)$.
%
%\begin{theorem}[Galerkin approximation]
%
%Assumptions (A), (B), (E) in force.
%Then $p^m \rightarrow p_X$ when
%$m \rightarrow \infty$
%
%\end{theorem}
%
%PROOF:
%
%
%We shall interpret the derivation of the projected laws $p^m$ as
%an extension of Galerkin-Faedo method to equations over
%differetiable manifolds, see \cite{Temam}.
%
\section{Acknowledgments}
The first named author wishes to thank
Bernard Hanzon and Fran\c{c}ois Le Gland
for their suggestions and their past work on the subject.


\begin{thebibliography}{99}

\bibitem{AMR82} Amari, S.-I. (1982). Differential geometry of Curved
Exponential Families-Curvature and Information Loss.{\em The Annals
of Statistics} {\bf 18}(2), 357-385.

\bibitem{AMR85} Amari, S. I. (1985). {\em Differential-Geometrical
Methods in Statistics}, Lecture Notes in Statistics 28,
Springer-Verlag, Berlin.

\bibitem{BafiPist} R. Bafico, G. Pistone, G-convergence of generators
and weak convergence of diffusions {\em Annales de l'institut Henri
Poincar\'e}, {\bf 1}, (1985), pagg. 1-13.

\bibitem{BRN78} Barndorff--Nielsen, O. E. (1978). {\em Information
and Exponential Families in Statistical Theory}. Wyley, New York.

\bibitem{BRC94} Barndorff--Nielsen O. E., and Cox D. R. (1994). {\em
Inference ans asymptotics}, Chapman \& Hall London

\bibitem{BRC89} Barndorff--Nielsen O. E., and Cox D. R. (1989). {\em
Asymptotic Techniques for use in Statistics.} Chapman \& Hall, London

\bibitem{BRN89} Barndorff--Nielsen, O. E. and Jupp, P. E. (1989).
Approximating Exponential Models. {\em Ann. Inst. Statist. Math.}
{\bf 41} 2, 247--267.

\bibitem{Brigo2} D. Brigo, On the nice behaviour of the Gaussian
projection filter with small observation noise,
Systems and Control Letters, vol. 26 (1995) 363--370.

\bibitem{BrHaLe} D. Brigo, B. Hanzon, F. Le Gland, A differential
geometric approach to nonlinear filtering: the projection filter,
{\em Publication Interne IRISA} {\bf 914}, IRISA, Rennes, 1995
(available at URL:~
ftp://ftp.irisa.fr/techreports/1995/PI-914.ps.Z). Reduced version
in: {\em Proceedings of the Conference on Decision and Control} (New
Orleans, 1995). Submitted to {\em IEEE Transactions on Automatic
Control}.

\bibitem{BrHaLe2} D. Brigo, B. Hanzon, F. Le Gland, On the
relationship between assumed density filters and projection filters,
TI Discussion Paper 7-96-18, 1996, Tinbergen Institute,
Amsterdam.

\bibitem{Brigo3} D. Brigo, New developments on the Gaussian projection
filter with small observation noise,
TI Discussion Paper 7-96-23, 1996, Tinbergen Institute, Amsterdam.

\bibitem{Volat} D. Brigo, B. Hanzon, Estimation of Stochastic
Volatility from bilateral exchange rates using projection filters,
working paper presented at the fourth Workshop of the {\em European
Network on System Identification} (ERNSI), Padua, June 7--9, 1995.

\bibitem{Brigo4} D. Brigo,  On  diffusions
with prescribed diffusion coefficients whose
densities evolve in  prescribed exponential families,
Internal Report CNR LADSEB, 02/96, March 1996.


\bibitem{EKT74} Ekeland, I. and Temam, R. (1974). {\em Analyse
convexe et probl\`{e}mes variationnels}. Dunod, Gauthier--Villars,
Paris.

\bibitem{El} K. D. Elworthy, {\em Stochastic Differential Equations
on Manifolds}, Cambridge University Press, Cambridge, 1982.

\bibitem{FrieA} A. Friedman, {\em Stochastic Differential Equations
and Applications}, vol. I (Academic Press, New York, 1975).

\bibitem{KRR58} Krasnosel'skii, M. A. and Rutickii, Ya. B. (1958).
{\em Convex Functions and Orlicz Spaces}. Fizmatgiz, Moskva. (In
Russian.) English translation (1961) Noordhoff, Groningen.

\bibitem{LNG95} Lang, S. (1995) {\em Differential
and Riemannian
Manifolds}. Springer-verlag, New York.

%\bibitem{MDS79} Madsen, L. T. (1979). {\em The geometry of
%statistical model - a generalization of curvature}, Research report
%79-1, Stat. Res. Unit, Danish Medical Res. Council.

\bibitem{MRR93} Murray, M. K. and Rice, J. W. (1993). {\em
Differential Geometry and Statistics}, Monographs in Statistics and
Applied Probability. {\bf 48} Chapman \& Hall.

\bibitem{Pardoux} E. Pardoux, Stochastic partial differential
equations and filtering of diffusion processes, {\em Stochastics},
{\bf 2} (1979), pp. 127- 167.

\bibitem{PSR94} Pistone G., and Rogantin M. P. (1994) {\em Geometria
differenziale: strumenti dell'inferenza statistica.} Corso di
Formazione, XXXVII Riunione Scientifica Societ\`a Italiana di
Statistica, Sanremo 9 aprile 1994.

\bibitem{PSR94a} Pistone G., and Rogantin M. P. (1994) {\em The
Transformation of the Non-Parametric Statistical Manifold under
Conditioning and Sampling.} 57th IMS Annual Meeting and 3rd World
Congress of the Bernoulli Society, Chapel Hill NC, June 20-25, 1994.

\bibitem{PSR96} Pistone G., and Rogantin, M. P. (1996)
{\em The Exponential Statistical Manifold: Mean
Parameters, Orthogonality and Space Transformation.} Rapporti
Interni del Dipartimento di Matematica del Politecnico di Torino No.
6/96.
Submitted to {\em Bernoulli}.

\bibitem{PSS92} Pistone, G., and Sempi, C. (1995) An
Infinite Dimensional Geometric Structure On the space of All the
Probability Measures Equivalent to a Given one. {\em The
Annals of Statistics} {\bf 23}(5) 1995.

\bibitem{RAR91} Rao, M. M., and Ren Z. D. (1991) {\em Theory of
Orlicz Spaces}. Dekker, New York.

\bibitem{RGN94} Rogantin M. P. (1994) {\em Modellizzazione
geometrica della dipendenza markoviana} Riunione scientifica del
gruppo ``Inferenza statistica: basi probabilistiche e sviluppi
metodologici'', Udine 16 -- 17 settembre 1994. Revised and submitted
to {\em Metron}, 1995

\bibitem{StroVarh} D.W.Stroock, S.R.S. Varadhan, (1979) {\em
Multidimensional diffusion processes} Springer Verlag, New York.

%\bibitem{Temam} Temam

\end{thebibliography}
\end{document}